\pdfoutput=1
\documentclass[draft=false, leqno]{article}
\usepackage{amsmath,amsthm,mathrsfs}
\usepackage[normal]{boilerart}

\newcommand{\BM}{\textit{BM}}

\theoremstyle{plain}
\newtheorem{red}[equation]{Reduction}

\title{A comment on the vanishing of rational motivic Borel--Moore homology}

\author{Clark Barwick \and Denis Nardin}

\date{}

\begin{document}

\maketitle

\begin{abstract} This note concerns a weak form of Parshin's conjecture, which states that the rational motivic Borel--Moore homology of a quasiprojective variety of dimension $m$ over a finite field in bidegree $(s,t)$ vanishes for $s>m+t$. It is shown that this conjecture holds if and only if the cyclic action on the motivic cohomology of an Artin--Schreier field extension in bidegree $(i,j)$ is trivial if $i<j$.
\end{abstract}


Let $k$ be a finite field of characteristic $p$; let $V$ be a quasiprojective variety of dimension $m$ over $k$. The conjecture of Beilinson--Parshin states that if $V$ is smooth and projective, then $K_i(V)\otimes\QQ=0$ for $i>0$; equivalently, the rational motivic cohomology $H^i(V,\QQ(j))$ vanishes unless $i=2j$. Equivalently, the conjecture states that for $V$ smooth and projective, $H^{\BM}_s(V,\QQ(t))$ vanishes unless $s=2t$.

We are interested in the following conjecture for arbitrary (i.e., not necessarily smooth or projective) $V$, which identifies a more restricted vanishing range:
\begin{cnj}\label{mainthm}
The rational motivic Borel--Moore homology $H^{\BM}_s(V,\QQ(t))$ vanishes if $s>m+t$.
\end{cnj}
\noindent Combined with usual vanishing results in motivic cohomology \cite[Th. 3.6 and Th. 19.3]{MR2242284}, this would imply that when $V$ is smooth (but not necessarily projective), one has (with $i=2m-s$ and $j=m-t$)
\[H^i(V,\QQ(j))=0\text{\; unless\; }i\in\left[j,j+m\right]\cap\left[j,2j\right].\]

Here is a conjecture concerning fields. Let $K$ be a perfect field of characteristic $p$, and let $L\coloneq K[y]/(y^p-y-a)$ be an Artin--Schreier extension, on which the cyclic group $C_p$ acts via $\goesto{y}{y+1}$.
\begin{cnj}\label{mainconj} The induced action of $C_p$ on $H^i(L,\QQ(j))$ is trivial for every $i<j$.
\end{cnj}
\noindent This would imply that $H^i(L,\QQ(j))$ vanishes in this range, so we may regard this as a kind of `ascent' property for motivic cohomology along Artin--Schreier covers.

The purpose of this note is to prove:
\begin{thm} \Cref{mainconj} implies \Cref{mainthm}
\end{thm}
\noindent The proof is an induction argument that reduces \Cref{mainthm} to \Cref{mainconj}. We are grateful to Joseph Ayoub, who kindly informed us that our previous formulation of this result was too strong.

\begin{nul} If $m=0$, Conjecture \ref{mainthm} (and indeed the Beilinson--Parshin Conjecture itself) follows from Quillen's computation of the $K$-theory of finite fields. When $m=1$, it follows from the celebrated computations of Harder. For the purpose of induction, we now assume this statement for quasiprojective varieties of dimension $<m$.
\end{nul}

\begin{nul} Choose an open immersion $\incto{V}{\overline{V}}$ into a projective variety of dimension $m$ such that the complement $\overline{V}-V$ (with its reduced scheme structure) is quasiprojective of positive codimension. The localization sequence
\[\cdots\to H^{\BM}_s(\overline{V}-V,\QQ(t))\to H^{\BM}_s(\overline{V},\QQ(t))\to H^{\BM}_s(V,\QQ(t))\to\cdots\]
now permits us to reduce to the case in which $V$ is projective. It suffices also to assume that $V$ is irreducible.
\end{nul}

Now we deploy the following result of Kiran Kedlaya:
\begin{thm}[Kedlaya, \protect{\cite[Theorem 1]{MR2092132}}]\label{thm:Kedlaya} Suppose $X$ a projective variety, pure of dimension $m$ over our finite field $k$. Suppose $L$ an ample line bundle on $X$, $D$ a closed subscheme of dimension less than $m$, and $S$ a $0$-dimensional subscheme of the regular locus not meeting $D$.

Then there exists a positive integer $r$ and an $(m+1)$-tuple of linearly independent sections of $L^{\otimes r}$ with no common zero such that the induced finite morphism
\[f\colon\fromto{X}{\PP H^0(X,L^{\otimes r})\cong\PP^m}\]
of $k$-schemes enjoys the following conditions.
\begin{enumerate}[(\ref{thm:Kedlaya}.1)]
\item If $\PP^{m-1}\cong H\subset\PP^m$ denotes the hyperplane at infinity, then $f$ is étale away from $H$.
\item The image $f(D)$ is contained in $H$.
\item The image $f(S)$ does not meet $H$.
\end{enumerate}
\end{thm}

\begin{nul} We thus obtain a finite morphism $f\colon\fromto{V}{\PP^m}$ that is étale over $\AA^m$. Let's write $Z\coloneq f^{-1}(H)$ and $U\coloneq f^{-1}(\AA^m)$; of course the latter is smooth.

The localization sequence
\[\cdots\to H^{\BM}_s(Z,\QQ(t))\to H^{\BM}_s(V,\QQ(t))\to H^{\BM}_s(U,\QQ(t))\to\cdots,\]
when combined with our induction hypothesis, reduces the problem to showing that the rational motivic cohomology
\[H^{i}(U,\QQ(j))\cong H^{\BM}_{2m-i}(U,\QQ(m-j))\]
vanishes whenever $i<j$.
\end{nul}

\begin{nul}\label{nul:passtofurthercover} At any stage, it will suffice to assume $U$ is connected, and moreover we will be free to pass to a further étale cover of $U$: indeed, if $g\colon\fromto{U'}{U}$ is a finite étale map, then the composite $g_\ast g^\ast\colon\fromto{H^i(U,\ZZ(j))}{H^i(U,\ZZ(j))}$ is multiplication by its degree. Hence
\[g^\ast\colon\fromto{H^i(U,\QQ(j))}{H^i(U',\QQ(j))}\]
is injective, and so it suffices to show that $H^i(U',\QQ(j))=0$ for $i<j$.
\end{nul}

\begin{nul} As a first application of \ref{nul:passtofurthercover}, if $f\colon\fromto{U}{\AA^m}$ is not Galois, we may pass to its Galois closure.

Harbater and van der Put show \cite[Example 5.3]{MR1928369} that a group is a finite quotient of the étale fundamental group of $\AA^m_{\overline{k}}$ (for $\overline{k}$ an algebraic closure of $k$) just in case it is a quasi-$p$-group. Hence by a second application of \ref{nul:passtofurthercover}, we may pass to a finite extension of $k$ and to connected components if necessary and thereby assume that $U$ is geometrically integral, and the Galois group $G$ of the Galois cover $f$ is a quasi-$p$-group.

By a third application of \ref{nul:passtofurthercover}, we may also pass to a finite extension of $k$ to ensure that the fiber of $f\colon\fromto{U}{\AA^m}$ over $0$ contains a rational point.
\end{nul}

\begin{nul} Since rational motivic cohomology satisfies étale descent, we have a convergent spectral sequence
\[E_2^{u,v}\cong H^u(G,H^v(U,\QQ(j)))\Rightarrow H^{u+v}(\AA^m_k,\QQ(j))\cong\begin{cases} \QQ&\text{if }u+v=0\text{\ and\ }j=0;\\ 0&\text{otherwise,}\end{cases}\]
by homotopy invariance and Quillen. Since $E_2^{u,v}$ vanishes unless $u=0$, we deduce that $H^{i}(U,\QQ(j))^G=0$ unless $i=j=0$.
\end{nul}

\begin{nul} The claim now is that $H^i(U,\QQ(j))=0$ is trivial when $i < j$; this is clearly true when
$G$ is the trivial group. Since $G$ is generated by elements of order a power of $p$ it suffices to show that every such element acts trivially. In particular, the conjecture will follow if for every Galois cover $U\to X$ of order $p^n$, the action of the Galois group on $H^i(U,\QQ(j))$ is trivial. We want to show that it suffices to check the case where $n=1$. We will prove this by induction on $n\ge2$.

Suppose we knew the above statement for Galois covers of order $p$, and let $g$ be a generator of the Galois group of $U$ over $X$. Suppose $n\ge2$. Then we can find $0<e<n$, so that both $e$ and $n-e$ are less than $n$. In particular, our thesis is true for $g^{p^e}$, that is the action of $g^{p^e}$ on $H^i(U,\QQ(j))$ is trivial. But then
\[H^i(U/g^{p^e},\QQ(j))=H^i(U,\QQ(j))^{g^{p^e}}=H^i(U,\QQ(j))\,.\]
Moreover, $g$ descends to an automorphism of $U/g^{p^e}$ of order $p^e$. Hence by our inductive hypothesis $g$ acts trivially on $H^i(U/g^{p^e},\QQ(j))=H^i(U,\QQ(j))$.

Since (as is well-known) Galois extensions of order $p$ are Artin–Schreier extensions, we may
now reduce to the following situation.

We suppose $A$ a smooth $k$-algebra, and we suppose that $A \subset B$ is an
Artin–-Schreier extension, so that $B\cong A[y]/(y^p-y-a)$. We assume that $T=\Spec A$ and $U=\Spec B$ are geometrically integral. Hence we may consider the subring $k[a]\subseteq A$; we note that since $U$ and $T$ are assumed geometrically integral, it follows that $a$ is not algebraic over $k$. Consequently, the function $a$ is a dominant, finite type morphism $a\colon\fromto{T}{\AA^1_k}$, and we have a pullback square
\begin{equation*}
\begin{tikzpicture}[baseline]
\matrix(m)[matrix of math nodes,
row sep=4ex, column sep=4ex,
text height=1.5ex, text depth=0.25ex]
{U & S \\
T & \AA^1_k, \\ };
\path[>=stealth,->,font=\scriptsize]
(m-1-1) edge node[above]{$b$} (m-1-2)
edge node[left]{$r$} (m-2-1)
(m-1-2) edge node[right]{$q$} (m-2-2)
(m-2-1) edge node[below]{$a$} (m-2-2);
\end{tikzpicture}
\end{equation*}
in which $S=\Spec k[x,y]/(y^p-y-x)$, and $q$ is the Artin--Schreier cover given by the inclusion $k[x]\subset k[x,y]/(y^p-y-x)$. (Of course $S\cong\AA^1_k$.)
\end{nul}

This, then, is our first reduction of \Cref{mainthm}:
\begin{red}\label{reduction1} The action of $C_p$ on $H^{i}(U,\QQ(j))$ is trivial if $i<j$.
\end{red}

\begin{nul} We now reduce the question to one of suitable function fields. That is, we claim that our induction hypothesis implies that if $V$ is smooth and geometrically irreducible, then $H^i(V,\QQ(j))\cong H^i(k(V),\QQ(j))$ for $i<j$. Indeed, for any nonempty open subset $W\subsetneqq V$, one has the localization sequence
\[
\to H_{2m-i}^{\BM}(V-W,\QQ(m-j))\to H^{i}(V,\QQ(j))\to H^{i}(W,\QQ(j))\to H_{2m-i-1}^{\BM}(V-W,\QQ(m-j))\to
\]
Let $c$ denote the codimension of $W$; note that $c\geq 1$, so that if $i<j$ then $2m-i-1>m-c+m-j$, whence by the induction hypothesis on the dimension,
\[
H_{2m-i}^{\BM}(V-W,\QQ(m-j))=H_{2m-i-1}^{\BM}(V-W,\QQ(m-j))=0.
\]
Consequently, one has an isomorphism
\[
H^{i}(V,\QQ(j))\cong H^{i}(W,\QQ(j))
\]
in this range. Passing to the colimit, one has $H^{i}(V,\QQ(j))\cong H^{i}(k(V),\QQ(j))$.
\end{nul}

\begin{red}\label{reduction2} The action of $C_p$ on $H^{i}(k(U),\QQ(j))$ is trivial if $i<j$.
\end{red}

\begin{nul} If $B$ is smooth over a perfect field $k$, then one may compare rational motivic cohomology of $B$ in the sense of Voevodsky with the Ext groups in the $\infty$-category $\categ{DM}(B;\QQ)$ of rational motives:
\[
H^i(B,\QQ(j))\cong[1_B,1_B(j)[i]]_{\categ{DM}(B;\QQ)}.
\]

In our case, we are interested in the situation in which $B$ is $\Spec$ of the function fields $k(T)$ and $k(U)$. We note that these fields are not perfect, but for any field $K$ with perfection $K^{\textit{perf}}$, the $\infty$-category $\categ{DM}(K;\QQ)$ is equivalent to $\categ{DM}({K^{\textit{perf}}};\QQ)$, so we are free to pass to the context originally contemplated by Voevodsky.

Consequently, we write $K\coloneq k(T)^{\textit{perf}}$, and $L\coloneq K(y)/(y^p-y-a)$.
\end{nul}

The task is thus to analyze the Galois action of the cyclic group $C_p$ on the rational motivic cohomology of $L\cong K[y]/(y^p-y-a)$ induced by the action $\goesto{y}{y+1}$. The final reduction of \Cref{mainthm} now is
\begin{red}\label{reduction3} The action of $C_p$ on $H^{i}(L,\QQ(j))$ is trivial if $i<j$.
\end{red}
\noindent This is \Cref{mainconj}. Equivalently, if we abuse notation slightly and write $L$ again for the Artin motive of $K\subset L$, then we have shown that \Cref{mainthm} would follow from the triviality of the action of $C_p$ on the cohomology $H^i(K,L(j))$ of the Artin--Tate motive $L(j)$ for $i<j$.


\DeclareFieldFormat{labelnumberwidth}{{#1\adddot\midsentence}}
\printbibliography
\addcontentsline{toc}{section}{References}

\end{document}